%% file: main.tex
\documentclass[12pt]{article}
\usepackage[utf8]{inputenc}

\input{micro.tex}

\usepackage{amsmath,booktabs}
\usepackage{listings}

\usepackage[vlined,linesnumbered,ruled]{algorithm2e}
\let\oldnl\nl
\newcommand{\nonl}{\renewcommand{\nl}{\let\nl\oldnl}}
\makeatother
\usepackage{mathrsfs}
\usepackage {eucal}
\usepackage{ amssymb }
\newcommand\scalemath[2]{\scalebox{#1}{\mbox{\ensuremath{\displaystyle #2}}}}

\title{Symmetry Lie Algebras of Varieties with Applications to Algebraic Statistics}
\author{Aida Maraj\thanks{Max Planck Institute of Molecular Cell Biology and Genetics, Center for Systems Biology Dresden, Technical University Dresden
  ({maraj@mpi-cbg.de}). AM was partially supported by the US National Science Foundation grant DMS-2306672.}  and Arpan Pal\thanks{University of Idaho (apal@uidaho.edu)}}

\date{}

\begin{document}

\maketitle

\begin{abstract}The motivation for this paper is to detect when an irreducible projective variety $V$ is not toric. This is achieved by analyzing a Lie group and a Lie algebra associated with $V$. If the dimension of~$V$ is strictly less than the dimension of the aforementioned objects, then $V$ is not a toric variety. We provide an algorithm to compute the Lie algebra of an irreducible variety and use it to present examples of non-toric statistical models in algebraic statistics.\end{abstract}


\section{Introduction}

The  paper is motivated  from the need to classify statistical models with toric structure in algebraic statistics, as many statistical models can be described as zero sets of polynomial ideals intersected with the probability space.  We say that a statistical model has a toric structure if its vanishing ideal is prime and generated by binomials, possibly after a linear change of variables. The toric structure in a statistical model is of interest due to the importance in applications: generating sets of toric ideals produce Markov bases and contribute in hypothesis testing algorithms~\cite{diaconis1998algebraic,Pistone.etal.2001}, facilitate maximum likelihood degree computations \cite{amendola2019maximum,boege2021reciprocal}, toric varieties are linked to smoothness criteria in exponential families~\cite{Geiger.etal.2001}, 
and the polytope associated to a toric model is useful when studying the existence of  maximum likelihood estimates~\cite{fienberg2012maximum}. 
Numerous papers, including \cite{beerenwinkel2007conjunctive,coons2021symmetrically,duarte2020equations,garcia2005algebraic,gorgen2022staged,hollering2021discrete,misra2021gaussian,sturmfels2005toric,sturmfels2019brownian}, bear witness to the interest in statistical models with toric structures.
While all these papers provide sufficient conditions for a statistical model to be toric,  the first statistical model  with non-toric structure was only recently provided by 
Nicklasson in~\cite{nicklasson2022toric}. 
This present paper proposes the symmetry Lie group associated to a homogeneous prime ideal  $I\subseteq  \mathbb C[x_1,\ldots,x_n]$ (\Cref{defn:sym_Lie_group}) and the symmetry Lie group of an irreducible projective variety $V\subseteq \C^n$ (\Cref{defn:sym_Lie_group_variety}) as  efficient ways to distinguish non-toric varieties. The two groups sit naturally in $\GL_n(\C)$ as stabilizers of the ideal/variety under  natural actions. They agree when $I$ is the vanishing ideal $I(V)$ of variety $V$. A fundamental observation of the paper is the following:
\begin{theorem}\label{thm:main} (see also \Cref{thm:main2})
Let $V$ be an irreducible projective variety with vanishing ideal~$I=I(V)$. 
Let  $G_{I}$  be the symmetry Lie group for $I$ as defined in \Cref{defn:sym_Lie_group}. If $\dim(G_I)< \dim(V)$ then~$V$ is not a toric variety.
\end{theorem}

Unfortunately, Lie groups tend to be challenging to compute. Since \Cref{thm:main} requires only the dimension of the symmetry Lie group, we are free to work with the Lie algebras of these groups. According to standard Lie theory literature, the dimension of a Lie group as a manifold is equal to the vector space dimension of its Lie algebra. Moreover, Lie algebras offer a friendlier structure, leading to the formulation of the following theorem and the subsequent development of an algorithm, which we implement using SageMath.

\begin{theorem}\label{thm:matrices} (\Cref{thm:matrices2} restated)
    Let $I\subseteq \mathbb C[x_1,\ldots, x_n]$ be a homogeneous prime ideal generated by polynomials of degree at most $d$.  Let  $ \mathscr{B}([I]_d)=\{f_1,\ldots, f_k\}$ be a finite  basis for the d-th graded component $[I]_d$ of $I$. Take $g\in M_n(\C)$ to be the $n\times n$ matrix with unknown entries $g_{ij}$.  For each $f_i\in \mathscr{B}([I]_d)$ consider the matrix 
\begin{align*}
    M_{i}(g):= \begin{pmatrix} \overrightarrow{f_1} &\overrightarrow{f_2} &\ldots & \overrightarrow{f_k} & \overrightarrow{g* f_i}
    \end{pmatrix},
\end{align*}
where the $*$ action is defined in \Cref{defn:*}, and $\overrightarrow{f_i}$ is the vector representation of polynomial~$f_i$ in $[R]_d$.
Then the symmetry Lie algebra for $I$ is the set of all matrices $g\in M_n(\C)$  such that $\mathrm{rank}{(M_i(g))}=k$ for $i=1,\ldots,k $. 
\end{theorem}

The observation and its implementation allow us to provide other statistical models with non-toric structure, shown in 
\Cref{sec:applications}. This includes  disproving  Conjecture~6.8 in~\cite{gorgen2022staged} which states that all staged tree models with one stage have toric structure. We also  provide an example of a Gaussian graphical model whose variety is not toric, which to the knowledge of the authors, is the first such example. We also discuss ways to use symmetry Lie algebras to find statistical models with toric structure. \\

\textbf{The structure of the paper}. 
In \Cref{sec:preliminaries} we recall definitions of a  toric variety, ideal with toric structure, graded components of a homogeneous polynomial ideal, Lie groups and Lie algebras. We also include relevant results on them. \Cref{sec:lie_groups} concerns introducing the symmetry Lie group of an ideal and of a variety. We show in  \Cref{prop:I=V} that the two definitions agree when working with the vanishing ideal of a variety. \Cref{sec:lie_groups}  concludes with a proof of \Cref{thm:main}. \Cref{sec:lie_algebras} concerns the Lie algebra for the symmetry Lie group of a homogeneous prime ideal, with a focus on describing this object practically via a group action given in  \Cref{defn:*}. Then we state \Cref{thm:main2}, which is the Lie algebra version of \Cref{thm:main}. The rest of \Cref{sec:lie_algebras} is in the service to facilitate computations of the symmetry Lie algebras of homogeneous prime ideals, hence proving \Cref{thm:matrices}. An algorithm for it is attached to this paper.  In \Cref{sec:applications} we apply the methods developed in this  paper to varieties arriving from staged tree models and Gaussian graphical models in algebraic statistics. We end with a discussion on other possible applications of symmetry Lie algebras.


\section{Preliminaries}\label{sec:preliminaries}
We start with  a list of notations that is kept uniform throughout the paper, unless specifically stated otherwise: 
\begin{itemize}
\item $\mb T_n$ is the $n$-dimensional algebraic torus isomorphic to $(\C\setminus \{0\})^n$
\item $R$ denotes the polynomial ring $\C[x_1,\ldots, x_n]$
\item  $M_n(\C)$ is the ring of $n\times n$ matrices with entries in $\C$
\item $\GL_n(\C)$ is the general linear group of invertible $n\times n$ matrices with entries in $\C$
\item $V$ denotes an irreducible  variety in $\C^n$
\item $I$ is a homogeneous prime ideal in $R$
\end{itemize}

We briefly recall the definitions of a toric variety, toric ideal, graded components of a polynomial ideal over a standard graded polynomial ring, Lie groups and Lie algebras. For details on toric varieties, ideals and graded components we recommend \cite{cox2013ideals,cox2011toric}. Books  \cite{fulton2013representation} and \cite{hall2013lie} offer comprehensive insights into Lie groups and their Lie algebras.  \\

\textbf{Toric varieties.}
Let $\mathcal{A}=\{\mathbf{a}_1,\ldots, \mathbf{a}_n\}$ be finite subset of the character lattice of the torus $\mathbb T_r$. Consider the map 
\begin{align}\label{eq:toric_variety}
    \Phi_{\mathcal{A}}: \mathbb T_r \longrightarrow (\mathbb C\setminus \{0\})^n, \quad \Phi_{\mathcal{A}}(\mathbf{t})=(\chi^{\mathbf{a}_1}(\mathbf{t}),\ldots,\chi^{\mathbf{a}_n}(\mathbf{t}))
\end{align}
 where $\chi^{\mathbf{a}_i}$ are morphisms from $\mathbb T_r$ to $\mathbb C\setminus \{0\}$ (characters) of the torus $\mathbb T_r$. If $\mathbb T_r=(\mathbb C\setminus \{0\})^r$, as in our situation, then each character has the form $\chi^{\mathbf{a}_i}(\mathbf t)=t_1^{a_{1i}}\cdots t_r^{a_{ri}}$ for $\mathcal{A}\subseteq \mathbb Z^r$. 
An affine variety $V$ is a {toric affine variety} if it is the Zariski closure of the image of $ \Phi_{\mathcal{A}}$, for some character lattice $\mathcal{A}$ of some torus $\mathbb T_r$.
In particular, one has that in a toric variety, the map in  \Cref{eq:toric_variety} induces the group action 
\begin{align}\label{eq:GA_toric_variety}
    \Tilde{\Phi}_{\mathcal{A}}: \mathbb T_r \times \mathbb C^n \longrightarrow \mathbb C^n, \quad \Tilde{\Phi}_{\mathcal{A}}(\mathbf{t},v)=(\chi^{\mathbf{a}_1}(\mathbf{t})v_1,\ldots,\chi^{\mathbf{a}_n}(\mathbf{t})v_n)
\end{align}
such that $\Tilde{\Phi}_{\mathcal{A}}(\mathbb T_r\times V)\subseteq V$. A common fact that we will use in this paper is that for toric variety $V$, the dimension of  the torus is equal or bigger to the dimension of the variety.\\

\textbf{Ideals with toric structure.} {{
A polynomial in $R$ is a  binomial if it is of the form $\mathbf x^{\mathbf{a}}-\mathbf x^\mathbf{b}$. 
An ideal $I\subseteq R$ is said to be \emph{toric} if it is prime and it has a generating set made of binomials. Equivalently, the ideal $I$ is toric if and only if it can be expressed as the kernel of a monomial map from a polynomial ring to a Laurent ring. The ideal $I$ has a \emph{binomial structure} if it has a generating set of binomials in variables $x_1,\ldots, x_n$, or if there exists an invertible linear change of variables such that the ideal has a binomial generating set in the new variables. 
 We say that $I$ has a \emph{toric structure} if it is prime and has a binomial structure.}}

\begin{example}\label{ex:coin_flip} (biased coin model in algebraic statistics, see \cite{gorgen2022staged}) 
The prime ideal $I=\langle x_1x_3-x_2x_3-x_2^2\rangle $ is clearly not a binomial ideal in variables $x_1,x_2,x_3$. Consider the change of variables $x_1=y_1$, $x_2=y_2+y_1$ and $x_3=y_3-y_2$. The generator of $I$ takes the form 
\[x_1x_3-x_2x_3-x_2^2= y_1(y_3-y_2)-(y_2-y_1)(y_3-y_2)-(y_2-y_1)^2= y_1^2-y_2y_3. \]
So, $I$ is toric in  $\mathbb C[y_1,y_2,y_3]$ under the presented linear change of variables. 
\end{example}
A common way to show that an irreducible variety is toric is by proving that it is the vanishing of a prime ideal with binomial structure. The end of this section lists relevant literature.\\ 

\begin{remark}
A binomial sometimes is defined as any expression of the form \( \mathbf{x}^{\mathbf{a}} - \lambda \mathbf{x}^{\mathbf{b}} \), where \( \lambda \) is a scalar. A unital binomial is a specific type of binomial where \( \lambda \) is restricted to either $0$ or $1$. Observe that a  binomial 
\begin{align}\label{eq:binomial}
   \mathbf{x}^{\mathbf{a}}-\lambda \mathbf{x}^\mathbf{b} \qquad \text{where } \lambda \neq 0,1
\end{align}can always be expressed into a unital binomial times a constant, under a simple linear change of variables. 
\end{remark}

\textbf{Graded components of an ideal.}
As a standard graded polynomial ring,  $R$ is a direct sum of its graded components; that is,  $R = \oplus _{d \geq 0} [R]_d$, where  \[[R]_d=\{p\in R \mid p \text{ is homogeneous polynomial of degree } d\}\cup \{0\}\] is the vector space of all  homogeneous polynomials in $R$ of degree $d$. We refer to the set of monomials of degree $d$ in $R$ as \emph{the standard basis} for $[R]_d$ and denote it $\mathscr{B}([R]_d)$. 
Similarly, a homogeneous ideal $I\subseteq R$ has $I = \oplus_{d \ge 0} [I]_d$, where the vector space \[[I]_d=\{p\in I \mid p \text{ is homogeneous polynomial of degree } d\}\cup \{0\}\] 
is  the \emph{d-th graded component of $I$}. The following proposition describes a basis for~$[I]_d$. 
\begin{proposition}\label{prop:basis}
Let $I\subseteq R$ be an ideal generated by homogeneous polynomials $p_1,\ldots,p_k$ of degrees $d_1,\ldots,d_k$, respectively. Then, for each $d\in \mathbb N$, \begin{align*}
\mathscr{S}([I]_d) =  \bigcup_{i\in [k], d_i\leq d} ^k\{mp_i \mid   m\in \mathscr{B}([R]_{d-d_i})\}
\end{align*}
is a  spanning set for the vector space $[I]_d$.  Consequently, any linearly independent spanning set~$\mathscr{B}([I]_d)\subseteq \mathscr{S}([I]_d)$ is a basis for $[I]_d$. 
\end{proposition}

\begin{proof}
An element $m_ip_i\in \mathscr{S}([I]_d)$ is a degree-$d$ polynomial in $I$, implying that $m_ip_i\in [I]_d$. As~$[I]_d$ is a linear space, linear combinations of these elements are in $[I]_d$, so $ \mathrm{Span}(\mathscr{S}([I]_d))\subseteq [I]_d$. Conversely, consider $p\in [I]_d\subseteq I$.  We can express  $p=\sum\limits_{i=1}^kf_ip_i$, where $f_i\in [R]_{d-d_i}$.
So,  $p$ is a linear combination of polynomials $m_ip_i$ for~$m_i\in \mathscr{B}([R]_{d-d_i})$, as desired.
\end{proof}

\textbf{Lie groups and Lie algebras.}   
A \emph{Lie group} is a group that is also a finite-dimensional smooth manifold, in which the group operations of multiplication and inversion are smooth maps. We will concern ourselves with matrix Lie groups that are complex manifolds.  The group  $\GL_n(\C)$ is a classical example of a matrix Lie group. {{We are interested in Lie subgroups of $\GL_n(\C)$.  Cartan's theorem states that any closed subgroup in $\GL_n(\C)$, referred in \cite{hall2013lie} as a matrix Lie group, is a Lie group. }}

The \emph{Lie algebra} $\g$ of a Lie group $G$ is the tangent space of $G$ at the identity.   For a Lie group $G\subseteq \GL_n(\C)$, its Lie algebra has the particular form (see \cite[Section~3.3]{hall2013lie})
\begin{align}\label{eq:lie_algebra}
\mathfrak{g}= \{g\in M_{n}(\C) \mid  e^{tg}\in G \text{ for all } t\in \mathbb R \}. 
\end{align}
{{The Lie algebra $\g$ is called a \emph{complex Lie algebra} if $ig \in \g$ for any $g\in \g$.
In this case $G$ is called a \emph{complex Lie group}. 
We are interested in the dimensions of these objects. 
}}
\begin{theorem}\label{thm:cartan}\cite[Corollary~3.45]{hall2013lie}
Let $G$ be a closed subgroup of $\GL_n(\C)$ with Lie algebra $\g$ and let $k$ be the dimension of $\g$ as a real vector space. Then $G$ is Lie group,{{ whose dimension as a real manifold is $k$.  If $\g$ is a \emph{complex Lie algebra}, then the dimension of $\g$ as a complex vector space is equal to the dimension of $G$ as a complex manifold, equal to $k/2$. }}
\end{theorem}
 {{The last line of the theorem comes from the fact that when $\g$ is complex Lie algebra, its dimension as real vector space is twice the dimension of $\g$ as a complex vector space. On the other side, treated as manifolds, the analogous statement holds for $G$.}}

A \emph{finite-dimensional Lie representation} $\Pi$ of a Lie group $G$ is a continuous group homomorphism from $G$ to $GL(V)$ where $V$ is a finite-dimensional complex vector space with~$\dim(V)\geq 1$.
\begin{theorem}\label{thm:representations}
\cite[Proposition 4.4]{hall2013lie} Let \( G \) be a matrix Lie group with Lie algebra \( \mathfrak{g} \) and let \( \Pi \) be a (finite-dimensional real or complex) representation of \( G \), acting on the space \( V \). Then,  there is a unique representation \( \pi \) of \( \mathfrak{g} \) acting on the same space such that
\(
\Pi(e^X) = e^{\pi(X)}
\)
for all \( X \in \mathfrak{g} \). The representation \( \pi \) can be computed as
\[
\pi(X) = \frac{d}{dt} \bigg|_{t=0} \Pi(e^{tX}).
\]
\end{theorem}
Representations are related to group actions; a representation $\Pi$ determines the group action of $G$ in $V$ by $g\cdot v=\Pi(g)v$, and a natural operation of the corresponding Lie algebra $\g$ in $V$, which we will later use. The following is a useful corollary that relates representations and actions. The result is found in many books; for details, we recommend \cite[Section 13.2]{humphreys2012linear}.
\begin{corollary}\label{cor:stable subspaces}
{{Let $G$ and $\g$ be as in \Cref{thm:representations}. Then $G$ and $\g$ keep stable the same subspaces of $V$.}}
\end{corollary}

\textbf{Relevant literature.}
\begin{itemize}
\item  Lie groups on varieties are not new to algebraic statistics. Draisma, Kuhnt and Zwiernik \cite{draisma2013groups,draisma2017automorphism} use the symmetry Lie group with respect to conjugation on varieties of symmetric matrices arriving from Gaussian graphical  models to advance problems related to their maximum likelihood estimate. 
\item Katthän, Micha{\l}ek and Miller in \cite{katthan2019polynomial}   address when an ideal is binomial under an automorphism. The set up is more general, not only for \emph{linear} change of variables and their test for binomiality after an automorphism is based on
comprehensive Gr\"obner bases (see Algorithm 4.2 and Remark 6.2 and Problem 3 in \cite{katthan2019polynomial}). Hence implementation is feasible only on simple examples. Our algorithm relies only on linear algebra computations.
\item Nicklasson's  approach in~\cite{nicklasson2022toric} uses matrix representation of quadratic forms to prove non-toricness. They represent a homogeneous quadratic generator $f$ as $xSx^T$, where $S$ is an $n\times n$ symmetric coefficient  matrix and $x$ is the vector of variables. If $f$ is a binomial then $S$ must have rank at most $4$. This rank does not change when $f$ is binomial under a linear transformation. 
\item Kreuzer and Walsh  provide an algorithm  in \cite{kreuzer2023computing} to compute the subideal of an ideal $I$ generated by all monomials and binomials in $I$  by computing their cellular decomposition. 
\item Kahle and Vill \cite{kahle2024efficiently} recently follow up on the approach in the present paper and give an algorithm which detects whether an ideal is toric under a linear change of variables; see \Cref{rem:algorithm}.
\end{itemize}


\section{Symmetry Lie Groups}\label{sec:lie_groups}
\textbf{Symmetry Lie groups of homogeneous prime ideals.}  Consider the group action of  $\GL_n(\C)$ on the polynomial ring $R=\C[x_1,\ldots,x_n]$, 
\begin{align}\label{eq:cdot_action}
\GL_n(\C)\times R \rightarrow R, \ (g,p(x))\mapsto g\cdot p(x)
\end{align}
given by $g\cdot p(x) = p(g^{-1}x)$.

Alternatively, one can think of this group action as  
substituting each variable $x_i$ in  $p(x)$ with $g_i x$, where $g_i$ is the $i$-th row of $g^{-1}$ and  $x=[x_1,\ldots,x_n]^{T}$ is the vector of variables. 

\begin{example} \label{ex:liegroup}
Take polynomial $p(x)=x_1^2+x_2^2+x_1x_2\in \C[x_1,x_2]$ and matrix $g^{-1}\coloneq \begin{bmatrix} g_{11} & g_{12} \\ g_{21} &g_{22} 
    \end{bmatrix}$. Then,
\[g\cdot p(x)=  (g_{11}^2+g_{21}^2+g_{11}g_{21})x_1^2+ (g_{12}^2+g_{22}^2+g_{12}g_{22})x_2^2+  (2g_{11}g_{12}+2g_{21}g_{22}+g_{11}g_{22}+g_{12}g_{21})x_1x_2. 
\]
\end{example}
Now we consider the group that acts as stabilizer for an ideal $I$ in $R$.
\begin{definition}\label{defn:sym_Lie_group} The {stabilizer of the ideal  $I\subseteq R$} is  \[G_I = \{g\in \GL_n(\C) \mid g\cdot p(x)\in I, \forall p(x) \in I\}.\]
\end{definition}
We will prove in \Cref{thm:sequence_Lie}  that when $I$ is a homogeneous prime ideal, $G_{I}$ is  a Lie subgroup of $\GL_n(\C)$, which we refer to as \emph{the symmetry Lie group} of $I$.   First, we simplify our problem and show that a generating set of $I$ is sufficient to determine the stabilizer $G_I$. 
 
 \begin{lemma}\label{lm:gens_Liegroup}
 Let $I=\langle p_1,\ldots,p_k\rangle\subseteq R$.  Then $G_I = \{g\in \GL_n(\C) \mid g\cdot p_i\in I, 1\leq i\leq k\}$.
 \end{lemma}
 
\begin{proof}
Denote $G'_I= \{g\in \GL_n(\C) \mid g\cdot p_i\in I, 1\leq i\leq k\}$. It is clear that $G_I\subseteq G'_I$. Take a  matrix $g\in G'_I$. A polynomial $p\in I$ has form $p=\sum\limits_{i=1}^kq_ip_i$ for some polynomials $q_1,\ldots, q_k\in R$.  Since  each $g\cdot p_i\in  I$, one has that  $g\cdot p=\sum\limits_{i=1}^k(g\cdot q_i)(g\cdot p_i)\in I$, and hence $g\in G_I$. 
 \end{proof}

\begin{example}\label{ex:lie group manifold} Consider the ideal $I=\langle p(x)\rangle $ where $p(x)=x_1^2+x_2^2+x_1x_2$ as in \Cref{ex:liegroup}.  Then $G_I$ is a $2$-dimensional manifold made of all invertible matrices $g\in \GL_n(\C)$ such that 
\begin{align*}
g_{11}^2+g_{21}^2+g_{11}g_{21}=g_{12}^2+g_{22}^2+g_{12}g_{22}=2g_{11}g_{12}+2g_{21}g_{22}+g_{11}g_{22}+g_{12}g_{21}.
\end{align*}
\end{example}

\textbf{Passing to $[I]_d$}. Note that action (\ref{eq:cdot_action}) is degree preserving; that is, the polynomials $p(x)$ and~$g\cdot p(x)$ have the same degree. This property allows us to safely define the group action on each graded component $[R]_d$, which can be treated as vector space with basis, say monomials of degree $d$ in $R$.  We can similarly talk about the stabilizer for vector space $[I]_d\subseteq [R]_d$.

\begin{definition}\label{defn:sym_Lie_group_d}
 Let~$I\subseteq R$ be a homogeneous ideal.  The stabilizer of the vector space $[I]_d$ is
  \begin{align*}
    G_{[I]_d} = \{g\in \GL_n(\C) \mid g\cdot p(x)\in [I]_d, \forall p(x) \in [I]_d\}. 
\end{align*}
\end{definition}

An analogous version to \Cref{lm:gens_Liegroup}  holds for $G_{[I]_d}$.

 \begin{lemma}\label{lm:gens_Liegroup2}
 Let $I$ be a homogeneous ideal in $R$ generated by the homogeneous polynomials $p_1,\ldots,p_k$ of degrees $d_1,\ldots, d_k$, respectively. For each $d\in \mathbb N$, let $\mathscr{B}([I]_d)$ be as in \Cref{prop:basis}.  Then, 
   \[G_{[I]_d} = \{g\in \GL_n(\C) \mid g\cdot p(x)\in [I]_d,  \forall p(x)\in \mathscr{B}([I]_d)\}.\]
 \end{lemma}

 Now let us look at the  example of the  maximal ideal in $R$. 

\begin{example}\label{ex:max_ideal}
Let $I=\langle x_1,\ldots, x_n\rangle$ in $R$. Then $[I]_d=[R]_d$ for all $d\in \N$.  So, 
\[G_I = \{g\in \GL_n(\C) \mid g\cdot x_i\in I, 1\leq i\leq n\}= \GL_n(\C),\]
and, for each $d\in \N$, we have
\(G_{[I]_d} = \{g\in \GL_n(\C) \mid g\cdot m\in [R]_d, \forall m\in \mathscr{B}([R]_d) \}= \GL_n(\C).\)
\end{example}

In \Cref{ex:max_ideal}, $G_I=G_{[I]_d}$ for any $d\in\N$. A similar statement holds in a  general setting. 

\begin{theorem}\label{thm:sequence_Lie}
Let $I$ be a homogeneous prime ideal in $R$ generated by homogeneous polynomials $p_1,\ldots,p_k$ of degrees $d_1,\ldots, d_k$, respectively. Then, $(G_{[I]_d})_{d\in \mathbb N}$ is a non-increasing sequence of Lie groups with respect to inclusion. Moreover, $G_I=G_{[I]_d}$ for $d\geq \mathrm{max}\{d_1,\ldots, d_k\}$.  
\end{theorem}

\begin{proof}
We will use \Cref{thm:cartan} to show that each $G_{[I]_d}$ is a Lie group, i.e., we will show that $G_{[I]_d}$ is a closed subgroup of $\GL_n(\C)$.
Pick~$d\in \mathbb N$. The identity matrix $Id_{n}\in G_{[I]_d}$ since $Id_{n}\cdot p(x)=p(x)\in [I]_{d}$ for any $p(x)\in [I]_d$. For any $g,h\in G_{[I]_d}$ and $p(x)\in [I]_{d}$, $gh\cdot p(x)=p((gh)^{-1}x)=p(h^{-1}g^{-1}x)=h\cdot (g\cdot p(x))\in [I]_{d}$. 
For $g\in G_{[I]_d}$, the map $\tilde{g}: [R]_d\rightarrow [R]_d$ given by  $\tilde{g}(p(x))=g\cdot p(x)$ is an invertible linear transformation with inverse~$\tilde{g^{-1}}$. Now $\tilde{g}([I]_d)\subseteq [I]_d$. Given that $[I]_d$ is of  finite dimension,  $\tilde{g}([I]_d)= [I]_d$, and so $g^{-1}([I]_d)=[I]_d$, which implies that $g^{-1}\in G_{[I]_d}$. Lastly, take a sequence of matrices $g_{i}\in G_{[I]_d}$ converging to~$g$. Then $g\cdot p(x)=(\lim\limits_{i\rightarrow \infty}g_i)\cdot p(x)=\lim\limits_{i\rightarrow \infty}(g_i\cdot p(x))\in [I]_d$.  Hence,  $G_{[I]_d}$ is a closed subgroup of~$\GL_n(\C)$. 

Next, we will show that the sequence of the Lie groups  $(G_{[I]_d})_{d\in \mathbb N}$ is non-increasing. We can safely assume $I\neq \langle x_1,\ldots, x_n\rangle$. The case $I=\langle x_1,\ldots, x_n\rangle$ was proved true in \Cref{ex:max_ideal}. 
Since $I\neq \langle x_1,\ldots, x_n\rangle$, there is some variable $x_t\in [R]_1$ not in $I$.  
Take $g\in G_{[I]_d}$.  Consider the linear form  $g^{-1}\cdot x_t\in [R]_1 $.  The property  $g\cdot p(x)q(x)=(g\cdot p(x))(g\cdot q(x))$ implies that for
 $p(x)\in [I]_{d-1}$, \[g\cdot ((g^{-1}\cdot x_t)p(x))=(g\cdot (g^{-1}\cdot x_t))(g\cdot p(x))=x_t(g\cdot p(x))\in [I]_d\subseteq I.\] Since $I$ is prime and $x_t\notin I$, one must have that $g\cdot p(x)\in I$. Our group action preserves the degree of a polynomial, so $g\cdot p(x)\in [I]_{d-1}$, and consequentially $g\in G_{[I]_{d-1}}$, which concludes that $G_{[I]_d}\subseteq G_{[I]_{d-1}}$.

For the final part of the theorem it is enough to show that  $ G_{[I]_d}\subseteq G_I$ for $d= \mathrm{max}\{d_1,\ldots, d_k\}$. Take $g\in G_{[I]_d}$. By \Cref{lm:gens_Liegroup}, it is enough to show that $g\cdot p_i(x)\in I$,  for any $1\leq i\leq k$. If $p_i(x)\in [I]_d\subseteq I$, we are done.
Otherwise, recall that $q(x)p_i(x)\in [I]_d$ for any standard basis element $q(x)\in \mathscr{B}([R]_{d-d_i})$. In particular, this is true for $q(x)= (g^{-1}\cdot x_t)^{d-d_i}$. Thus we have 
\[g\cdot ((g^{-1}\cdot x_t)^{d-d_i}p_i(x))= x_t^{d-d_i}(g\cdot p_i(x))\in [I]_d\subseteq I,\]
where $x_t$ is one of the variables not in $I$, as earlier. Given that $I$ is  a prime ideal and $x_t\notin I$, one has that $g\cdot p_i(x)\in I$, as desired. In particular,  $G_I$ is a Lie group.
\end{proof}

\textbf{Symmetry Lie groups of varieties.} One can similarly define the symmetry Lie group of an irreducible variety. Consider  $\GL_n(\C)$ acting on $\C^n$ with the rule:  \[\text{for each  $v\in \C^n$ and $g=( g_{ij})_{n\times n}\in \GL_n(\C)$, one has $ g \bullet v = {g}v$},\]
where  the point $v$ is interpreted as an $n\times 1$ vector.

\begin{definition} \label{defn:sym_Lie_group_variety}
The {stabilizer of  variety $V\subseteq \C^n$} is 
\(
    G_V = \{g\in \GL_n(\C) \mid  g\bullet v\in V, \forall v \in V\}.
\)
\end{definition}
\noindent This definition is chosen to coincide with the stabilizer of the vanishing ideal of that~variety.

\begin{proposition}\label{prop:I=V}
Let  $V$ variety in $\mathbb C^{n}$ be affine variety with homogeneous vanishing ideal $I=I(V) \subseteq  R$. Then~$G_{I}=G_V$. In particular, $G_V$ is a Lie group when $V$ is irreducible. \end{proposition}

\begin{proof}
The definitions of the two group actions imply that for any $p\in I$ and $v\in V$, one has \begin{align}\label{eq:group_properties}
p(g \bullet v)=(g^{-1}\cdot p)(v).\end{align} Hence, 
\begin{align*}
g\in G_V  &\longleftrightarrow p(g \bullet v)=0 \text{ for all } p\in I,  v \in V  \quad \text{(by the definition of $G_V$)}\\
&\longleftrightarrow  (g^{-1}\cdot p)(v)=0 \text{ for all } p\in I, \ v\in V   \quad \text{(by \Cref{eq:group_properties})} \\
&\longleftrightarrow    g^{-1}\in G_I  \quad \text{(by the definition of $G_I$).}\\
&\longleftrightarrow    g\in G_I  \quad \text{(since  $G_I$ is a group)}.
\end{align*}
The final part is true from the correspondence between prime homogeneous ideals and irreducible projective varieties.
\end{proof}
We will refer to $G_V$ as \emph{the symmetry Lie group} of $V$. Now, we are ready to prove  \Cref{thm:main}.

\begin{proof}[Proof of \Cref{thm:main}]
Suppose that $V\subseteq \C^n$ is a toric variety with torus $\mb T_r$. So, $ \dim(V)= \dim(\mb T_r)$.
Next, we show that 
 $\dim (\mb T_r)\leq \dim(G_V)$ by providing an embedding of $\mb T_r$ in $G_V$. 
Start with the embedding~$\iota$ 
\begin{align*}
    \iota: \mb T_n\rightarrow \GL_n(\C), \ (t_1,t_2,\ldots, t_n)\xrightarrow{} \begin{pmatrix}
      {t_1} & 0 & \ldots & 0\\
      0 & {t_2}  & \ldots & 0\\
      \vdots & \vdots & \ddots & \vdots\\
      0 & 0 & \ldots & {t_n} 
    \end{pmatrix}
\end{align*}
and consider the composition $ \iota \circ {\Phi_{\mathcal{A}}}$ for some $\Phi_{\mathcal{A}}$ in \Cref{eq:toric_variety}. For $t\in \mb T_r$ and $v\in V$, 
\[\iota \circ  {\Phi_{\mathcal{A}}}(\mathbf{t}) \bullet v =  (t^{a_1}v_1, \cdots, t^{a_n}v_n)= \Tilde{\Phi}_{\mathcal{A}}(\mathbf{t},v)\in V, \text{ for } \Tilde{\Phi}_{\mathcal{A}} \text{ as in \Cref{eq:GA_toric_variety}} .\] Hence $\iota \circ  {\Phi_{\mathcal{A}}}(\mb T_r)\subseteq  G_V$, which makes  $ \iota \circ  {\Phi_{\mathcal{A}}}$  a desired embedding. 
The chain of inequalities $\dim(V)= \dim(\mb T_r)\leq \dim(G_V)$ concludes the proof. 
\end{proof}
  
The rest of the article concerns with only homogeneous prime ideals in polynomial rings. \Cref{prop:I=V} allows for interpretations of the results to irreducible varieties. 
  

\section{Symmetry Lie Algebras}\label{sec:lie_algebras}
 
Let $G_I$ be the symmetry Lie group of a homogeneous prime ideal $I$ in $R$. The \emph{Lie algebra} of $G_I$, which we will refer to it as the \emph{symmetry Lie algebra} of the ideal $I$ and denote by $\g_I$, is
\begin{align*}
\g_I = \{g\in M_n(\C)\mid e^{tg} \in G_I,  \forall t\in \mathbb R\}.
\end{align*} 
We will describe $\g_I$ as stabilizer of $I$ under the action induced by $\cdot$ on $R$ defined as follows:  
\begin{align}\label{eq:star}
  *  \colon M_{n}(\C) \times  R\rightarrow R, \quad g *  p(x)= \dfrac{d}{dt}|_{t=0}(e^{tg} \cdot p(x)).
\end{align}

\begin{theorem} \label{thm:Liealgebra I}
 Let  $I=\langle p_1(x),\ldots, p_k(x)\rangle$  be a homogeneous prime ideal and let \[d\geq \max\{\deg(p_1(x)), \ldots, \deg(p_k(x))\}.\] Then,  $\g_I=\g_{[I]_d}$. It is precisely  
 \begin{align}
     \g_{[I]_d}=  \{g\in M_n(\C)\mid g *  [I]_d\subseteq [I]_d \}. 
 \end{align}
\end{theorem}

\begin{proof}
$G_I=G_{[I]_d}$ by \Cref{thm:sequence_Lie}. Hence, we conclude  that $\g_I=\g_{[I]_d}$. 

We now focus on $[I]_d\subseteq [R]_d$ and use its vector space structure. The $\cdot$ action restricted to $[R]_d$ uniquely defines a rational finite complex Lie representation with rational entries \[
 \Pi\colon \GL_{n}(\C)\rightarrow \GL([R]_d).\]
By \Cref{thm:representations}, we obtain the representation on the respective Lie algebras
\begin{align*}
 \pi\colon M_{n}(\C)\rightarrow \mathfrak{gl}([R]_d), \quad \pi(g)= \dfrac{d}{dt}|_{t=0}\Pi(e^{tg}).
\end{align*}
whose associated action on $[R]_d$ by $M_n(\C)$ is precisely the $ * $ operation given in \Cref{eq:star} restricted to $[R]_d$.  
We can restrict the representations $\Pi$ and $\pi$ to $G_{[I]_d}$ and $\g_{[I]_d}$, respectively. 
By \Cref{cor:stable subspaces},~$G_I$ and $\g_I$ must fix the same subspaces of $[R]_d$. The definition of $G_{[I]_d}$ as the largest subgroup fixing $[I]_d$  gives that $\g_I$ must also be the largest Lie subalgebra fixing $[I]_d$.  Hence, $\g_{[I]_d}=\g_I$ must be the Lie algebra
\( \g_{[I]_d} = \{g\in M_n(\C)\mid g *  [I]_d\subseteq [I]_d \}, \)  as desired. 
\end{proof}

We observe that the $*$ action can be described  independently of the differential in (\ref{eq:star}). This approach makes computing $\g_I$ more feasible. 

\begin{proposition}\label{defn:*} The  $*$ action of $M_{n}(\C)$ on $R$ is given by the rules: 
    \begin{enumerate}
        \item $g * c = 0$  for any constant  $c\in R$,
        \item $g * x_i = -\sum_{j=1}^ng_{ij}\cdot x_j$  for any variable $x_i\in R$,
        \item $g*(p_1p_2) = (g*p_1)p_2+p_1(g*p_2)$,  for any $p_1,p_2\in R$,
    \end{enumerate}
 which extended linearly to $R$ fully determine it. 
\end{proposition}

\begin{proof}
Following the definition of $*$-action, it is enough to check that \Cref{defn:*} holds  for monomials in~$R$. 

First note that for $p(x)\in R$ and $g\in M_n(\C)$,  $e^{tg}\cdot p(x)=p(e^{-tg}x)$. 

Take $g\in M_n(\C)$. 
If~$p(x)=c$ then
  $g*c={\dfrac{d}{dt}|}_{t=0}(e^{tg}\cdot c)={\dfrac{d}{dt}|}_{t=0}(c)=0$. If $p=x_i$ then
${{\dfrac{d}{dt}}|}_{t=0}(e^{tg}\cdot x_i)={\dfrac{d}{dt}|}_{t=0} \sum\limits_{j=1}^n(e^{-tg})_{ij}x_j= \sum\limits_{j=1}^n-g_{ij}x_j$.

Suppose \Cref{defn:*}  holds for any monomial of degree at most $d$. A monomial of degree $d+1$ in $R$ has form $x_im$, for some variable $x_i$  monomial $m$  of degree $d$ in~$R$. We have,
\begin{align*}g*x_im={\dfrac{d}{dt}|}_{t=0}(e^{tg}\cdot x_im)
&={\dfrac{d}{dt}|}_{t=0}(e^{tg}\cdot x_i)(e^{tg}\cdot m)\quad  (\textit{property (3) of } \cdot \textit{action})  
\\&= x_i{\dfrac{d}{dt}|}_{t=0}(e^{tg}\cdot m) + m{\dfrac{d}{dt}|}_{t=0}(e^{tg}\cdot x_i)  \quad (\textit{diff. rule}) 
\\ &= x_i (g*m)+m(g*x_i) \quad  (\textit{induction hypothesis}) 
\end{align*}
\end{proof}

\Cref{defn:*} together with \Cref{lm:lie_algebra_gens} will be particularly useful in computing $\g_I$.
\begin{corollary}\label{lm:lie_algebra_gens}
Let $I\subseteq R$ be a homogeneous prime ideal. For $d$ such $G_I=G_{[I]_d}$, let $\mathscr{B}([I]_d)=\{f_1(x),\ldots,f_k(x)\}$ be a basis of the vector space $[I]_d$. Then 
\begin{align*}
     \g_{I}&= \{g\in M_n(\C)\mid g *  f_i(x)\in [I]_d,  i=1,\ldots,k\}. 
 \end{align*}
\end{corollary}

\begin{proof}
  Denote $\g'_{[I]_d}=\{g\in M_n(\C)\mid g *  f_i(x)\in [I]_d, i=1,\ldots,k\}$. Clearly, $\g_{[I]_d}\subseteq \g'_{[I]_d}$. Conversely,
take~$g\in\g'_{[I]_d}$. For arbitrary $p(x)\in [I]_d$ we have  $p(x)=\sum\limits_{i=1}^k c_if_i(x)$. So, $g *  p(x)=\sum\limits_{i=1}^k c_i(g *  f_i(x))\in [I]_d$ since each term of the sum is in $[I]_d$.  
\end{proof}

\begin{example}
Consider ideal $I$ in \Cref{ex:liegroup} and $g\in M_2(\mathbb C)$. Then
   \begin{align*}g*(x_1^2+&x_2^2+x_1x_2) =g*x_1^2+g*x_2^2+g*x_1x_2  \\ =&2x_1(-g_{11}x_1-g_{12}x_2)+ 2x_2(-g_{21}x_1-g_{22}x_2)   +(-g_{11}x_1-g_{12}x_{2})x_2+x_1(-g_{21}x_1-g_{22}x_{2})\\
   =&-(2g_{11} + g_{21})x_1^2 - (2g_{12} + 2g_{21} + g_{11} + g_{22})x_1x_2 - (2g_{22} + g_{12})x_2^2, \end{align*}
and  $\mathfrak{g}_I$ is the $2$-dimensional space \(\{g\in M_2(\mathbb C) \mid 2g_{11} + g_{21}=2g_{12} + 2g_{21} + g_{11} + g_{22}=2g_{22} + g_{12}\}.\) 
\end{example}

Now we are ready to rephrase \Cref{thm:main} in terms of symmetry Lie algebras.
\begin{theorem}\label{thm:main2}
Let $I$ be a homogeneous prime ideal and let $\g_I$ be its symmetry Lie algebra. If $\dim(\g_I)< \dim(I)$, then  $V(I)$ is not a toric variety. 
\end{theorem}

{{\begin{proof}
Given $g\in M_n(\C)$ and homogeneous polynomial $f\in R$. Note that the  expression $g*f(x)$ has coefficients that are homogeneous linear forms in the variables $g_{ij}$. Hence, when applying matrix $ig$ to~$f(x)$, one has $(ig)* f(x)=i(g*f(x))$. Using the definition of $\g_I$ in \Cref{lm:lie_algebra_gens}, one has $g\in \g_I$ if and only if $ig \in \g_I$. Hence, $\g_I$ is a complex Lie algebra.  Lastly, combining  \Cref{thm:cartan} and \Cref{thm:main} we complete the proof.
\end{proof}}}


For the remainder of this section, we will employ the $*$ action defined in \Cref{defn:*} to develop an algorithm for computing the symmetry Lie algebra of a homogeneous prime ideal. To accomplish this, we will revisit the graded components of an ideal and reinterpret polynomials in them as vectors in a linear space.

Fix an order on the monomials $\mathscr{B}([R]_d)$. The paper and the algorithm use the  reverse lexicographic order, but any order works. The \emph{vector representation} of a polynomial $p(x)\in  [R]_d$ with respect to this ordered basis, is the vector $\overrightarrow{p}\in \C^{\binom{n+d-1}{d}}$ of coefficients of $p(x)$ in the chosen order.  For instance, $p(x)=x_1^2+2x_1x_3-x_2x_3\in [\C[x_1,x_2,x_3]]_2$ has  $\overrightarrow{p}=[0 \ -1 \ 2 \ 0 \ 0 \ 1]^T$. 

We are ready for \Cref{thm:matrices}.  
\begin{theorem}\label{thm:matrices2}
    Let $I\subseteq R$ be a homogeneous prime ideal  generated by  polynomials  of degree at most $d$.  Let  $ \mathscr{B}([I]_d)=\{f_1,\ldots, f_k\}$ be a basis for $[I]_d$. Let  $g\in M_n(\C)$ be the $n\times n$ matrix whose entries  $g_{ij}$ are unknown.  For each $f_i\in \mathscr{B}([I]_d)$ consider the matrix 
\begin{align*}
    M_{i}(g):= \begin{pmatrix} \overrightarrow{f_1} &\overrightarrow{f_2} &\ldots & \overrightarrow{f_k} & \overrightarrow{g* f_i}
    \end{pmatrix}.
\end{align*}
Then,
\( \g_I=\{ g\in M_n(\C) \mid \mathrm{rank}{(M_i(g))}=k, \ \text{for } i=1,\ldots,k  \}. \)
\end{theorem}

\begin{proof}
For each $d\in \mathbb N$, let $\g_{[I]_d}$ be the Lie algebra of the symmetry Lie group $G_{[I]_d}$. 
By \Cref{thm:Liealgebra I}, $\g_I=\g_{[I]_d}$. We use \Cref{lm:lie_algebra_gens}. 
For each $1\leq i\leq k$, one has that $g*f_i\in [I]_d$ if and only if $g*f_i$ is a linear combination of the basis elements $\mathscr{B}([I]_d)$, if and only if $\overrightarrow{g*f_i}$ is a linear combination of~$\overrightarrow{f_1}, \ldots, \overrightarrow{f_k}$, if and only if the matrix $M_i(g)$ has rank $k$, as desired. 
\end{proof}
An implementation of the algorithm written in SageMath \cite{sagemath} can be found on GitHub at the following URL: \url{https://github.com/arpan-pal/Toric_via_symmetry}. The code uses the basis~$\mathscr{B}([R]_d)$ described in \Cref{prop:basis} in the reverse lexicographic order. To speed up the selection of a base $\mathscr{B}([I]_d)$, it is recommended to input a minimal  generating set  for the ideal~$I$, when possible.

\section{Applications to Algebraic Statistics}\label{sec:applications}

As discussed in the introduction, the classification of statistical models with toric structure is of interest in algebraic statistics. When the vanishing ideal of a statistical model is not toric, a preferred method to check if its variety is toric is by searching for linear transformations under which the vanishing ideal becomes toric. The binomial structure has been successfully investigated in phylogenetics \cite{sturmfels2005toric,sturmfels2019brownian}, staged tree models \cite{gorgen2022staged,nicklasson2022toric}, and several Bayesian networks \cite{beerenwinkel2007conjunctive,garcia2005algebraic,hollering2021discrete}. 
In this section, we apply \Cref{thm:main2} and its implementation in \Cref{thm:matrices} on the ideals for staged tree models and Gaussian graphical models. There is much more space left for exploration, and we encourage the interested reader to delve into it. \\

\textbf{Staged tree models.}  Staged tree models are discrete statistical models encoding relationships between events. 
They are realizable as rooted trees with colored vertices and labeled edges directed away from the root. Vertices represent events, edge labels represent conditional probabilities, and the colors on the vertices represent an equivalence relation. 
Vertices of the same color have the same outgoing edge labels. 
We use $\theta_{ij}$ to denote the label associated with an edge~$[i,j]$. A key constraint is that the sum of the labels of all edges emanating from the same vertex in a staged tree must be equal to one. The staged tree model is defined as the set of points in $\mathbb R^n$ parametrized by multiplying edge labels along the root-to-leaf paths $\lambda_1, \ldots, \lambda_n$ in the staged tree $\mathcal{T}$. In algebro-geometric terms, the staged tree model consists of points inside the toric variety $V(\ker\varphi_\mathcal{T})$, where
\begin{equation}\label{eq:varphi}
\varphi_{\mathcal{T}}: \R[x_1,\ldots, x_n] \to\R[\Theta,z]/\langle \sum\limits_{j}\theta_{ij}-z\rangle, \quad x_r \mapsto z^{n-\ell(\lambda_r)}\prod_{[i,j]\in E(\lambda_r)}\theta_{ij} \text{ for } r=1,\ldots,n.
\end{equation} 
For an introduction to staged tree models, we refer the reader to \cite{duarte2020equations}. Detailed information on toric staged tree models after a linear change of variables can be found in \cite{gorgen2022staged}. The latter paper poses several open questions, two of which we address here. The first example of a staged tree model that has been shown to be non-toric is credited to Nicklasson \cite{nicklasson2022toric}. The code associated with \cite{gorgen2022staged} includes an implementation of \Cref{eq:varphi} that we utilize for faster computations.
\begin{figure}[h]
\centering
\begin{minipage}[b]{0.45\linewidth}
\centering
  \includegraphics[scale=0.7]{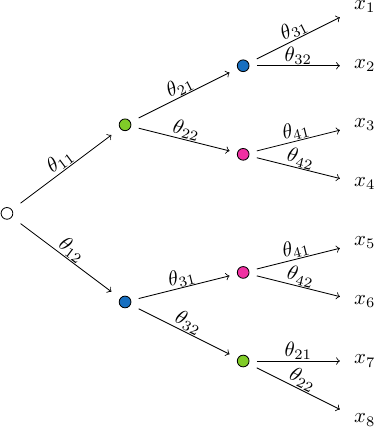} 
\caption{Staged tree from \Cref{ex:liealgebra_staged}}
\label{fig:staged tree 1}
\end{minipage}
\quad
\begin{minipage}[b]{0.45\linewidth}
\centering
   \includegraphics[scale=0.73]{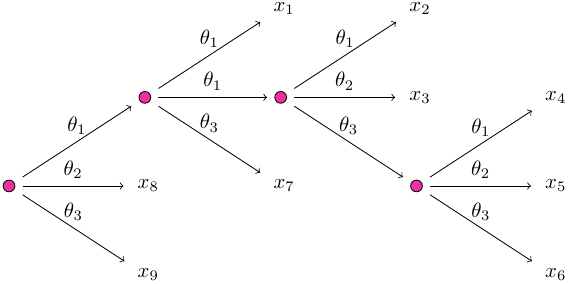}
\caption{Staged tree from \Cref{ex:caterpillar tree}}
\label{fig:staged tree 2}
\end{minipage}%
\end{figure}
\begin{example}\label{ex:liealgebra_staged}
The discussion section in \cite{gorgen2022staged} raises the question of whether $\ker(\varphi_{\mathcal T})$ 
for the staged tree, $\mathcal{T}$, in \Cref{fig:staged tree 1}, becomes toric after a linear change of variables. We provide a negative answer. The ideal  $\ker\varphi_{\mathcal{T}}\subseteq \C[x_1,\ldots,x_8]$ has dimension $5$  and generated by the $3$ quadratics \(p_1= (x_1+x_2)x_8-(x_3+x_4)x_7, \ p_2 =(x_7+x_8)x_1-(x_5+x_6)x_2,\text{ and } p_3 = x_3x_6-x_4x_5.\) So,  $\mathscr{B}([I]_2)=\{p_1,p_2,p_3\}$, and
{\small
 \begin{align*}
&  \overrightarrow{p_1}= [0 , 0   , 0,    0, 0   , 0   , 1   , 1   , 0   , 0   , 0   , -1   , -1   , 0   , 0   , 0   , 0   , 0   , 0   , 0   , 0   , 0   , 0   , 0   , 0   , 0   , 0   , 0   , 0   , 0   , 0   , 0   , 0   , 0   , 0   , 0 ]^{T},\\
 &  \overrightarrow{p_2} =  [0   , 0   , 0   , 0   , 0   , 0   , 0   , 1   , 0   , 0   , 0   , 0   , 0   , 0   , 1   , 0   , 0   , 0   , 0   , -1   , 0   , 0   , 0   , 0   , -1   , 0   , 0   , 0   , 0   , 0   , 0   , 0   , 0   , 0   , 0   , 0]^{T},\\
&  \overrightarrow{p_3}=  [0   , 0   , 0   , 0   , 0   , 0   , 0   , 0   , 0   , 0   , 0   , 0   , 0  \  0  \  0  \  0  \  0   , 0   , 1   , 0   , 0   , 0  , -1   , 0   , 0   , 0,  \  0   , 0   , 0   , 0   , 0   , 0   , 0   , 0   , 0   , 0]^{T}.\\
 &  \overrightarrow{g\star p_1}=  [ g_{18} + g_{28}  ,g_{17} + g_{27} - g_{38} - g_{48}   , g_{16} + g_{26}   , g_{15} + g_{25}   , g_{14} + g_{24} - g_{78}
   , g_{13} + g_{23} - g_{78},\\ & g_{12}  +  g_{22} + g_{88},    g_{11} + g_{21} + g_{88}  , -g_{37} - g_{47}   , -g_{36} - g_{46}   , -g_{35} - g_{45}   , -g_{34} - g_{44} - g_{77}, -g_{33} - \\ &g_{43} -  g_{77}   , -g_{32} - g_{42} + g_{87}   , -g_{31} - g_{41} + g_{87}   , 0   , 0   , -g_{76}   , -g_{76}   , g_{86}   , g_{86}   , 0   , -g_{75},  -g_{75},  g_{85},   g_{85} ,\\ & -g_{74},  -g_{73} - g_{74}   , -g_{72} + g_{84}   , -g_{71} + g_{84}   , -g_{73}   , -g_{72} + g_{83}   , -g_{71} + g_{83}   , g_{82}   , g_{81} + g_{82}   , g_{81}]^{T}.
\end{align*}  }
Compute $M_ 1(g)=\begin{pmatrix}   \overrightarrow{p_1} & \overrightarrow{p_2}&    \overrightarrow{p_3}&   \overrightarrow{g\star p_1}\end{pmatrix}$, and similarly  $M_ 2(g)$ and~$M_ 3(g)$.
Via Algorithm \ref{thm:matrices} we obtain that $\g_{\ker(\varphi_\mathcal{T})}$ is the $4$-dimensional vector space generated by
\[ \scalemath{0.7}{\begin{pmatrix}
1 & 0 & 0 & 0 & 0 & 0 & 0 & 0 \\
0 & 1 & 0 & 0 & 0 & 0 & 0 & 0 \\
0 & 0 & 0 & 0 & 0 & 0 & 0 & 0 \\
0 & 0 & 1 & 1 & 0 & 0 & 0 & 0 \\
0 & 0 & 0 & 0 & 0 & 0 & 0 & 0 \\
0 & 0 & 0 & 0 & 1 & 1 & 0 & 0 \\
0 & 0 & 0 & 0 & 0 & 0 & 1 & 0 \\
0 & 0 & 0 & 0 & 0 & 0 & 0 & 1
\end{pmatrix}}, \
\scalemath{0.65}{\begin{pmatrix}
0 & 0 & 0 & 0 & 0 & 0 & 0 & 0 \\
0 & 0 & 0 & 0 & 0 & 0 & 0 & 0 \\
0 & 0 & 1 & 0 & 0 & 0 & 0 & 0 \\
0 & 0 & -1 & 0 & 0 & 0 & 0 & 0 \\
0 & 0 & 0 & 0 & 0 & 0 & 0 & 0 \\
0 & 0 & 0 & 0 & -1 & -1 & 0 & 0 \\
0 & 0 & 0 & 0 & 0 & 0 & -1 & 0 \\
0 & 0 & 0 & 0 & 0 & 0 & 0 & -1
\end{pmatrix}}, \
\scalemath{0.65}{\begin{pmatrix}
0 & 0 & 0 & 0 & 0 & 0 & 0 & 0 \\
0 & 0 & 0 & 0 & 0 & 0 & 0 & 0 \\
0 & 0 & 0 & 1 & 0 & 0 & 0 & 0 \\
0 & 0 & 0 & -1 & 0 & 0 & 0 & 0 \\
0 & 0 & 0 & 0 & 0 & 1 & 0 & 0 \\
0 & 0 & 0 & 0 & 0 & -1 & 0 & 0 \\
0 & 0 & 0 & 0 & 0 & 0 & 0 & 0 \\
0 & 0 & 0 & 0 & 0 & 0 & 0 & 0
\end{pmatrix}}, \
\scalemath{0.7}{\begin{pmatrix}
0 & 0 & 0 & 0 & 0 & 0 & 0 & 0 \\
0 & 0 & 0 & 0 & 0 & 0 & 0 & 0 \\
0 & 0 & 0 & 0 & 0 & 0 & 0 & 0 \\
0 & 0 & 0 & 0 & 0 & 0 & 0 & 0 \\
0 & 0 & 0 & 0 & 1 & 0 & 0 & 0 \\
0 & 0 & 0 & 0 & 0 & 1 & 0 & 0 \\
0 & 0 & 0 & 0 & 0 & 0 & 1 & 0 \\
0 & 0 & 0 & 0 & 0 & 0 & 0 & 1
\end{pmatrix}}.
\]
By \Cref{thm:main2},   $\ker(\varphi_{\mathcal{T}})$ is not toric.
\end{example}

\begin{example} \label{ex:caterpillar tree}
 We use \Cref{thm:main} on the ideal of the staged tree, $\mathcal{T}$, in \Cref{fig:staged tree 2}, to disprove Conjecture~6.8 in \cite{gorgen2022staged} which states that under a linear transformation all staged tree models with one stage have toric vanishing ideal. 
 The ideal $\ker(\varphi_\mathcal{T})$ of the one-stage tree in \Cref{fig:staged tree 2} is of dimension $3$ and  minimally generated by  the two by two minors of the matrix 
 \[
\begin{pmatrix}
x_1+ . . . +x_6 & x_1 & x_2 & x_4 \\
x_8 & x_{2}+. . . +x_6 & x_3 & x_5 \\
x_9 & x_7 & x_4+x_5+x_6 & x_6
\end{pmatrix}. 
\] 
Using  \Cref{thm:matrices} and simplifications (details of implementation in GitHub)  we get that its symmetry Lie algebra is a $2$-dimensional vector space generated by 
 \[\scalemath{0.7}{\begin{pmatrix}
1 & 0 & 0 & 0 & 0 & 0 & 0 & 0 & 0 \\
0 & 1 & 0 & 0 & 0 & 0 & 0 & 0 & 0 \\
0 & 0 & 1 & 0 & 0 & 0 & 0 & 0 & 0 \\
0 & 0 & 0 & 1 & 0 & 0 & 0 & 0 & 0 \\
0 & 0 & 0 & 0 & 1 & 0 & 0 & 0 & 0 \\
0 & 0 & 0 & 0 & 0 & 1 & 0 & 0 & 0 \\
0 & 0 & 0 & 0 & 0 & 0 & 1 & 0 & 0 \\
0 & 0 & 0 & 0 & 0 & 0 & 0 & 1 & 0 \\
0 & 0 & 0 & 0 & 0 & 0 & 0 & 0 & 1
\end{pmatrix}}, \ \scalemath{0.7}{\begin{pmatrix}
-2 & 0 & 0 & 0 & 0 & 0 & 0 & 0 & 0 \\
0 & -3 & 0 & 0 & 0 & 0 & 0 & 0 & 0 \\
0 & 0 & -3 & 0 & 0 & 0 & 0 & 0 & 0 \\
0 & 1 & 0 & -4 & 0 & 0 & 0 & 0 & 0 \\
0 & 0 & 1 & 0 & -4 & 0 & 0 & 0 & 0 \\
0 & 0 & 0 & 2 & 2 & -2 & 0 & 0 & 0 \\
1 & 1 & 1 & 1 & 1 & 1 & -1 & 0 & 0 \\
0 & 0 & 0 & 0 & 0 & 0 & 0 & -1 & 0 \\
1 & 1 & 1 & 1 & 1 & 1 & 1 & 1 & 0
\end{pmatrix} }.\]
 By \Cref{thm:main}, $\ker(\varphi_\mathcal{T})$ is not toric.
\end{example}

A naïve hope is for a staged tree model whose tree contains a subtree with non-toric structure to not have a toric structure. As the example below shows,  this is not true.

\begin{remark}
Consider the  maximal one staged tree model of depth $4$ with each vertex of degree three. It contains the staged tree in \Cref{fig:staged tree 2} as its subtree. By  \cite[Lemma~6.1]{gorgen2022staged}, all maximal one stage trees have toric vanishing ideals, so  this ideal is toric -- it has a minimal generating set of $66$ linear binomials and $75$ quadratic binomials. 
\end{remark}

However, we suspect the following to hold. 
\begin{conjecture}
Let $\mathcal T$ be a staged tree with non-toric $V(\ker \varphi_{\mathcal T})$ that uses color set $S$. Let $\mathcal T'$ be a staged tree such that its restriction to color set $S$ gives $\mathcal T$. Then  $\ker \varphi_{\mathcal T'}$ is not toric under any linear change of variables.
\end{conjecture}
This is mainly based on the work done on gluing to leaves of a balanced trees, trees with the subtree inclusion property  (SIP) \cite{gorgen2022staged} and colors different from colors used by the balanced tree. The result there was positive; SIP trees are toric under a linear change of variables, and a proper combination of these linear transformations gave a linear transformation in all the tree, but it suggested that non-toricity of ideals of staged trees, will be preserved when such trees are glued in leaves of a balanced tree, as long as the latter one has colors different from the rest. \\

\textbf{Gaussian graphical models.}  A Gaussian graphical model is a collection of multivariate Gaussian distributions in which a graph encodes conditional independence relations among the random variables. Its set of concentration matrices, that is, the inverses of covariance matrices, is a linear space of symmetric matrices intersected with the cone of positive definite matrices. For a graph $G$, this linear space $\mathcal{L}_G$ is the set of all $n\times n$ symmetric matrices $K=(k_{ij})_{1\leq i,j\leq n}$ with $k_{ij}=0$ if $[i,j]$ is not an edge in $G$. The inverse of $K$, when this is invertible, is the covariance matrix $\Sigma=(\sigma_{ij})_{1\leq i,j\leq n}$ of the corresponding Gasussian distribution. The Zariski closure of all symmetric matrices $\Sigma\in \mathbb{R}^{n\times n}$ such that $\Sigma^{-1}\in \mathcal{L}_G$ is not linear, and most often  not a friendly variety. One can compute it  as the vanishing of the kernel of the rational map  \begin{align} \label{eq:minors}
\rho_G: \mathbb{R}[\Sigma] \to \mathbb{R}(K), \ 
\rho_G(\sigma_{ij})= \frac{(-1)^{i+j}K_{[n]\setminus \{i\},[n]\setminus\{j\}}}{\det (K)}, 
\end{align}
where~$K_{[n]\setminus \{i\},[n]\setminus\{j\}}$ is the~$ij$-th minor of the symmetric matrix~$K$.
\begin{figure}[h]
\centering
\begin{minipage}[b]{0.45\linewidth}
\centering
\scalemath{0.9}{\begin{tikzpicture}
\draw[very thick] (0,2) -- (0,0);
\draw[very thick] (2,0) -- (2,2);
\draw[very thick] (0,2) -- (2,2);
\draw[very thick] (0,0) -- (2,0);
	\draw [] plot [ mark=square*] coordinates {(0,2)};
		\draw [] plot [ mark=square*] coordinates {(0,0)};
	\draw [] plot [  mark=square*] coordinates { (2,0)};
	\draw [] plot [only marks,  mark = square*] coordinates {(2,2)};
	\node[left] at (0,2) {1};
	\node[left] at (0,0) {2};
	\node[right] at (2,0) {3};
	\node[right] at (2,2) {4};
\end{tikzpicture} }
\caption{Graphical model in \Cref{ex:gaussian}}
\label{fig:four cycle}
\end{minipage}
\qquad
\begin{minipage}[b]{0.45\linewidth}
\centering
 \scalemath{0.9}{ \begin{tikzpicture}
\draw[very thick, cyan](0,0) -- (2,1) -- (4,0);
\draw[very thick, magenta] (2,1) -- (4,0);
	\draw [green] plot [only marks, mark size = 3, mark=*] coordinates {(0,0)};
	\draw [blue] plot [only marks, mark size=3, mark=*] coordinates { (2,1)};
	\draw [green] plot [only marks, mark size = 3, mark=*] coordinates {(4,0)};
	\node[left] at (0,0) {1};
	\node[above] at (2,1) {2};
	\node[right] at (4,0) {3};
\end{tikzpicture}}
\caption{Graphical model in \Cref{ex:toric}}
\label{fig:colored graph}
\end{minipage}%
\end{figure} Misra and Sullivant show in \cite{misra2021gaussian} that block graphs produce toric Gaussian graphical models. In the next example we show that the ideal $\rho_\mathcal{G}$ arriving from a four cycle is not toric under any linear change of variables. To the knowledge of the authors this is the first example of a Gaussian graphical model that has been shown to be non-toric.

\begin{example} \label{ex:gaussian}
Consider the four cycle $G$ with edges $[1,2],[2,3],[3,4],[1,4]$ as in \Cref{fig:four cycle}. 
The ideal $\ker(\rho_G)$ of the Gaussian graphical model in \Cref{fig:four cycle} is of dimension $8$ and generated by \begin{align*}&p_1=\sigma_{23}\sigma_{14}\sigma_{24}-\sigma_{13}\sigma_{24}^2-\sigma_{22}\sigma_{14}\sigma_{34}+\sigma_{12}\sigma_{24}\sigma_{34}+\sigma_{22}\sigma_{13}\sigma_{44}-\sigma_{12}\sigma_{23}\sigma_{44},\\ & p_2=\sigma_{13}\sigma_{23}\sigma_{14}-\sigma_{24}\sigma_{13}^2-\sigma_{12}\sigma_{33}\sigma_{14}+\sigma_{11}\sigma_{33}\sigma_{24}+\sigma_{12}\sigma_{13}\sigma_{34}-\sigma_{11}\sigma_{23}\sigma_{34}.\end{align*}
Following the reverse lexicographic ordering on the indices of the variables; i.e.  $\sigma_{11}>\sigma_{12}>\sigma_{22}>\cdots > \sigma_{44}$,  the symmetry Lie algebra is a $4$-dimensional vector space generated by 
 \[\scalemath{0.55}{\begin{pmatrix}
1 & 0 & 0 & 0 & 0 & 0 & 0 & 0 & 0 & 0 \\
0 & 0 & 0 & 0 & 0 & 0 & 0 & 0 & 0 & 0 \\
0 & 0 & -1 & 0 & 0 & 0 & 0 & 0 & 0 & 0 \\
0 & 0 & 0 & 0 & 0 & 0 & 0 & 0 & 0 & 0 \\
0 & 0 & 0 & 0 & -1 & 0 & 0 & 0 & 0 & 0 \\
0 & 0 & 0 & 0 & 0 & -1 & 0 & 0 & 0 & 0 \\
0 & 0 & 0 & 0 & 0 & 0 & 0 & 0 & 0 & 0 \\
0 & 0 & 0 & 0 & 0 & 0 & 0 & -1 & 0 & 0 \\
0 & 0 & 0 & 0 & 0 & 0 & 0 & 0 & -1 & 0 \\
0 & 0 & 0 & 0 & 0 & 0 & 0 & 0 & 0 & -1
\end{pmatrix}}, \ \scalemath{0.55}{\begin{pmatrix}
0 & 0 & 0 & 0 & 0 & 0 & 0 & 0 & 0 & 0 \\
0 & 1 & 0 & 0 & 0 & 0 & 0 & 0 & 0 & 0 \\
0 & 0 & 2 & 0 & 0 & 0 & 0 & 0 & 0 & 0 \\
0 & 0 & 0 & 0 & 0 & 0 & 0 & 0 & 0 & 0 \\
0 & 0 & 0 & 0 & 1 & 0 & 0 & 0 & 0 & 0 \\
0 & 0 & 0 & 0 & 0 & 0 & 0 & 0 & 0 & 0 \\
0 & 0 & 0 & 0 & 0 & 0 & 0 & 0 & 0 & 0 \\
0 & 0 & 0 & 0 & 0 & 0 & 0 & 1 & 0 & 0 \\
0 & 0 & 0 & 0 & 0 & 0 & 0 & 0 & 0 & 0 \\
0 & 0 & 0 & 0 & 0 & 0 & 0 & 0 & 0 & 0
\end{pmatrix} },\  \scalemath{0.55}{\begin{pmatrix}
0 & 0 & 0 & 0 & 0 & 0 & 0 & 0 & 0 & 0 \\
0 & 0 & 0 & 0 & 0 & 0 & 0 & 0 & 0 & 0 \\
0 & 0 & 0 & 0 & 0 & 0 & 0 & 0 & 0 & 0 \\
0 & 0 & 0 & 1 & 0 & 0 & 0 & 0 & 0 & 0 \\
0 & 0 & 0 & 0 & 1 & 0 & 0 & 0 & 0 & 0 \\
0 & 0 & 0 & 0 & 0 & 2 & 0 & 0 & 0 & 0 \\
0 & 0 & 0 & 0 & 0 & 0 & 0 & 0 & 0 & 0 \\
0 & 0 & 0 & 0 & 0 & 0 & 0 & 0 & 0 & 0 \\
0 & 0 & 0 & 0 & 0 & 0 & 0 & 0 & 1 & 0 \\
0 & 0 & 0 & 0 & 0 & 0 & 0 & 0 & 0 & 0
\end{pmatrix}},\  \scalemath{0.55}{
  \begin{pmatrix}
0 & 0 & 0 & 0 & 0 & 0 & 0 & 0 & 0 & 0 \\
0 & 0 & 0 & 0 & 0 & 0 & 0 & 0 & 0 & 0 \\
0 & 0 & 0 & 0 & 0 & 0 & 0 & 0 & 0 & 0 \\
0 & 0 & 0 & 0 & 0 & 0 & 0 & 0 & 0 & 0 \\
0 & 0 & 0 & 0 & 0 & 0 & 0 & 0 & 0 & 0 \\
0 & 0 & 0 & 0 & 0 & 0 & 0 & 0 & 0 & 0 \\
0 & 0 & 0 & 0 & 0 & 0 & 1 & 0 & 0 & 0 \\
0 & 0 & 0 & 0 & 0 & 0 & 0 & 1 & 0 & 0 \\
0 & 0 & 0 & 0 & 0 & 0 & 0 & 0 & 1 & 0 \\
0 & 0 & 0 & 0 & 0 & 0 & 0 & 0 & 0 & 2
\end{pmatrix}}.\]
By \Cref{thm:main}, $\ker(\rho_G)$ is not toric under any linear change of variables.
\end{example}
A natural question emerges:
\begin{question}
 Let $G$ be a  graph, such that $V(\ker(\rho_G))$ is not toric. Let $G'$ be a  graph with~$G$  as subgraph. Is it true that $V(\ker(\rho_{G'}))$ is not a toric variety?    
\end{question}
\textbf{Colored Gaussian graphical models.} Colored  Gaussian graphical models are generalizations of Gaussian graphical models. The graph is colored, and in addition to $k_{ij}=0$ whenever $[i,j]$ is a missing edge in the graph, one has that $k_{ii}=k_{jj}$ when vertices $i$ and $j$ have the same color, and $k_{ij}=k_{uv}$ when the edges $[i,j]$ and $[u,v]$ have the same color. The vanishing ideal for a colored graph $\mathcal G$ is the kernel of the map in \Cref{eq:minors} adapted to the subset of parameters in $K$.  For an introduction to these models we recommend \cite{uhler2011} and for work on the colored Gaussian graphical models with toric vanishing ideals see \cite{coons2021symmetrically}.  Note that \Cref{ex:gaussian} also serves as the first example of a \emph{colored} Gaussian graphical model whose vanishing ideal is not toric under any linear change of variables. 
We use symmetry Lie algebras on an ideal arriving from a colored Gaussian graphical model to discover its toric structure. 

\begin{example}\label{ex:toric}
Consider ideal $I=\ker(\rho_\mathcal{G})\subseteq R=\C[\sigma_{11},\sigma_{12},\ldots, \sigma_{33}]$ for the colored graph $\mathcal{G}$  in \Cref{fig:colored graph}. This is a $4$-dimensional ideal  generated by 
\(p_1=\sigma_{12}\sigma_{13}-\sigma_{11}\sigma_{23} \text{ and } p_2= \sigma_{12}^2-\sigma_{11}\sigma_{22}-\sigma_{13}^2+\sigma_{11}\sigma_{33}.\) 
Following the reverse lexicographic ordering on the indices of the variables; i.e.  $\sigma_{11}>\sigma_{12}>\sigma_{22}>\cdots > \sigma_{44}$, the symmetry Lie algebra of~$I$ is $11$-dimensional vector space generated by 
\begin{align*} 
&\scalemath{0.6}{\begin{pmatrix} 1 & 0 & 0 & 0 & 0 & 0 \\
0 & 0 & 0 & 0 & 0 & 0 \\
0 & 0 & 0 & 0 & 0 & 0 \\
0 & 0 & 0 & 0 & 0 & 0 \\
0 & 0 & 0 & 0 & -1 & 0 \\
0 & 0 & 1 & 0 & 0 & -1
\end{pmatrix}}, \ 
\scalemath{0.6}{\begin{pmatrix}0 & 0 & 0 & 0 & 0 & 0 \\
1 & 0 & 0 & 0 & 0 & 0 \\
0 & 0 & 0 & 0 & 0 & 0 \\
0 & 0 & 0 & 0 & 0 & 0 \\
0 & 0 & 0 & 1 & 0 & 0 \\
0 & -2 & 0 & 0 & 0 & 0
\end{pmatrix}}, \ 
\scalemath{0.6}{\begin{pmatrix}0 & 0 & 0 & 0 & 0 & 0 \\
0 & 1 & 0 & 0 & 0 & 0 \\
0 & 0 & 0 & 0 & 0 & 0 \\
0 & 0 & 0 & 1 & 0 & 0 \\
0 & 0 & 0 & 0 & 2 & 0 \\
0 & 0 & -2 & 0 & 0 & 2
\end{pmatrix}}, \
 \scalemath{0.6}{\begin{pmatrix}0 & 0 & 0 & 0 & 0 & 0 \\
0 & 0 & 0 & 1 & 0 & 0 \\
0 & 0 & 0 & 0 & 0 & 0 \\
0 & -1 & 0 & 0 & 0 & 0 \\
0 & 0 & -1 & 0 & 0 & 1 \\
0 & 0 & 0 & 0 & -4 & 0
\end{pmatrix}},  \ 
\scalemath{0.6}{\begin{pmatrix} 0 & 0 & 0 & 0 & 0 & 0 \\
0 & 0 & 0 & 0 & 0 & 0 \\
0 & 0 & 0 & 0 & 0 & 0 \\
1 & 0 & 0 & 0 & 0 & 0 \\
0 & 1 & 0 & 0 & 0 & 0 \\
0 & 0 & 0 & 2 & 0 & 0
\end{pmatrix}}, \\
&\scalemath{0.6}{\begin{pmatrix} 0 & 0 & 0 & 0 & 0 & 0 \\
0 & 0 & 0 & 0 & 0 & 0 \\
1 & 0 & 0 & 0 & 0 & 0 \\
0 & 0 & 0 & 0 & 0 & 0 \\
0 & 0 & 0 & 0 & 0 & 0 \\
1 & 0 & 0 & 0 & 0 & 0
\end{pmatrix}}, \
 \scalemath{0.6}{\begin{pmatrix} 0 & 0 & 0 & 0 & 0 & 0 \\
0 & 0 & 0 & 0 & 0 & 0 \\
0 & 1 & 0 & 0 & 0 & 0 \\
0 & 0 & 0 & 0 & 0 & 0 \\
0 & 0 & 0 & 0 & 0 & 0 \\
0 & 1 & 0 & 0 & 0 & 0
\end{pmatrix}},\ 
 \scalemath{0.6}{\begin{pmatrix} 0 & 0 & 0 & 0 & 0 & 0 \\
0 & 0 & 0 & 0 & 0 & 0 \\
0 & 0 & 1 & 0 & 0 & 0 \\
0 & 0 & 0 & 0 & 0 & 0 \\
0 & 0 & 0 & 0 & 0 & 0 \\
0 & 0 & 1 & 0 & 0 & 0
\end{pmatrix}},\  
\scalemath{0.6}{\begin{pmatrix} 0 & 0 & 0 & 0 & 0 & 0 \\
0 & 0 & 0 & 0 & 0 & 0 \\
0 & 0 & 0 & 1 & 0 & 0 \\
0 & 0 & 0 & 0 & 0 & 0 \\
0 & 0 & 0 & 0 & 0 & 0 \\
0 & 0 & 0 & 1 & 0 & 0
\end{pmatrix}},\ 
 \scalemath{0.6}{\begin{pmatrix}  0 & 0 & 0 & 0 & 0 & 0 \\
0 & 0 & 0 & 0 & 0 & 0 \\
0 & 0 & 0 & 0 & 1 & 0 \\
0 & 0 & 0 & 0 & 0 & 0 \\
0 & 0 & 0 & 0 & 0 & 0 \\
0 & 0 & 0 & 0 & 1 & 0
\end{pmatrix}},\ 
 \scalemath{0.6}{\begin{pmatrix} 0 & 0 & 0 & 0 & 0 & 0 \\
0 & 0 & 0 & 0 & 0 & 0 \\
0 & 0 & 0 & 0 & 0 & 1 \\
0 & 0 & 0 & 0 & 0 & 0 \\
0 & 0 & 0 & 0 & 0 & 0 \\
0 & 0 & 0 & 0 & 0 & 1
\end{pmatrix}}.
\end{align*}
Since $\dim(\g_{I})\geq \dim(I)$, \Cref{thm:main} is inconclusive. We will show that $I$ is toric after an appropriate linear change of variables by further analyzing~$\g_I$.
We are in search for a torus in~$\g_I$ of dimension at least $4$.   
Consider the invertible matrix
\begin{align*}
  B= \scalemath{0.8}{\left(\begin{array}{rrrrrr}
0 & 0 & 0 & 0 & 0 & 1 \\
0 & -i & i & 0 & 0 & 0 \\
0 & 0 & 0 & 0 & 1 & 0 \\
0 & 1 & 1 & 0 & 0 & 0 \\
1 & 0 & 0 & 1 & 0 & 0 \\
2 i & 0 & 0 & -2 i & 1 & 0
\end{array}\right)}. 
\end{align*}
Apply a change of basis in $\g_I$ with respect to $B$. The new basis (so $B^{-1}AB$ for each basis element~$A$) is listed below, and it contains  exactly $4$ diagonal matrices, revealing the embedded $4$-dimensional torus: 
\begin{align*}&\scalemath{0.6}{\left(\begin{array}{rrrrrr}
-1 & 0 & 0 & 0 & 0 & 0 \\
0 & 0 & 0 & 0 & 0 & 0 \\
0 & 0 & 0 & 0 & 0 & 0 \\
0 & 0 & 0 & -1 & 0 & 0 \\
0 & 0 & 0 & 0 & 0 & 0 \\
0 & 0 & 0 & 0 & 0 & 1
\end{array}\right)},\scalemath{0.6}{
\left(\begin{array}{rrrrrr}
0 & 1 & 0 & 0 & 0 & 0 \\
0 & 0 & 0 & 0 & 0 & \frac{1}{2} i \\
0 & 0 & 0 & 0 & 0 & -\frac{1}{2} i \\
0 & 0 & 1 & 0 & 0 & 0 \\
0 & 0 & 0 & 0 & 0 & 0 \\
0 & 0 & 0 & 0 & 0 & 0
\end{array}\right)},\scalemath{0.6}{
\left(\begin{array}{rrrrrr}
2 & 0 & 0 & 0 & 0 & 0 \\
0 & 1 & 0 & 0 & 0 & 0 \\
0 & 0 & 1 & 0 & 0 & 0 \\
0 & 0 & 0 & 2 & 0 & 0 \\
0 & 0 & 0 & 0 & 0 & 0 \\
0 & 0 & 0 & 0 & 0 & 0
\end{array}\right)},\scalemath{0.6}{
\left(\begin{array}{rrrrrr}
2 i & 0 & 0 & 0 & 0 & 0 \\
0 & i & 0 & 0 & 0 & 0 \\
0 & 0 & -i & 0 & 0 & 0 \\
0 & 0 & 0 & -2 i & 0 & 0 \\
0 & 0 & 0 & 0 & 0 & 0 \\
0 & 0 & 0 & 0 & 0 & 0
\end{array}\right)}, \scalemath{0.6}{
\left(\begin{array}{rrrrrr}
0 & 0 & 0 & 0 & 0 & 0 \\
0 & 0 & 0 & 0 & 0 & 0 \\
0 & 0 & 0 & 0 & 0 & 0 \\
0 & 0 & 0 & 0 & 0 & 0 \\
0 & 0 & 0 & 0 & 0 & 1 \\
0 & 0 & 0 & 0 & 0 & 0
\end{array}\right),}\\
&
\scalemath{0.6}{\left(\begin{array}{rrrrrr}
0 & 0 & 0 & 0 & 0 & 0 \\
0 & 0 & 0 & 0 & 0 & 0 \\
0 & 0 & 0 & 0 & 0 & 0 \\
0 & 0 & 0 & 0 & 0 & 0 \\
0 & -i & i & 0 & 0 & 0 \\
0 & 0 & 0 & 0 & 0 & 0
\end{array}\right)},\scalemath{0.6}{
\left(\begin{array}{rrrrrr}
0 & 0 & 0 & 0 & 0 & 0 \\
0 & 0 & 0 & 0 & 0 & 0 \\
0 & 0 & 0 & 0 & 0 & 0 \\
0 & 0 & 0 & 0 & 0 & 0 \\
0 & 0 & 0 & 0 & 1 & 0 \\
0 & 0 & 0 & 0 & 0 & 0
\end{array}\right)}, \scalemath{0.6}{
\left(\begin{array}{rrrrrr}
0 & 0 & 0 & 0 & 0 & 0 \\
0 & 0 & 0 & 0 & 0 & 0 \\
0 & 0 & 0 & 0 & 0 & 0 \\
0 & 0 & 0 & 0 & 0 & 0 \\
0 & 1 & 1 & 0 & 0 & 0 \\
0 & 0 & 0 & 0 & 0 & 0
\end{array}\right)},\scalemath{0.6}{
\left(\begin{array}{rrrrrr}
0 & 0 & 0 & 0 & 0 & 0 \\
0 & 0 & 0 & 0 & 0 & 0 \\
0 & 0 & 0 & 0 & 0 & 0 \\
0 & 0 & 0 & 0 & 0 & 0 \\
1 & 0 & 0 & 1 & 0 & 0 \\
0 & 0 & 0 & 0 & 0 & 0
\end{array}\right)},\scalemath{0.55}{
\left(\begin{array}{rrrrrr}
0 & 0 & 0 & 0 & 0 & 0 \\
0 & 0 & 0 & 0 & 0 & 0 \\
0 & 0 & 0 & 0 & 0 & 0 \\
0 & 0 & 0 & 0 & 0 & 0 \\
2 i & 0 & 0 & -2 i & 1 & 0 \\
0 & 0 & 0 & 0 & 0 & 0
\end{array}\right)},\scalemath{0.55}{
\left(\begin{array}{rrrrrr}
0 & -i & 0 & 0 & 0 & 0 \\
0 & 0 & 0 & 0 & 0 & \frac{1}{2} \\
0 & 0 & 0 & 0 & 0 & \frac{1}{2} \\
0 & 0 & i & 0 & 0 & 0 \\
0 & 0 & 0 & 0 & 0 & 0 \\
0 & 0 & 0 & 0 & 0 & 0 
\end{array}\right)}.\end{align*}
Consider the change of variables in $R$ induced by rows of $B$: 
\begin{align*} & \sigma_{11} \mapsto s_{33} & 
& \sigma_{12} \mapsto -is_{12}+is_{22} &
& \sigma_{22} \mapsto s_{23}\\
& \sigma_{13} \mapsto s_{12}+s_{22} &
& \sigma_{23} \mapsto s_{11}+s_{13} &
& \sigma_{33} \mapsto 2is_{11}-2is_{13}+s_{23}.
\end{align*}
This change of variables sends the original generating polynomials $p_1,p_2$ of $I$ to 
\[p_1' = -is_{12}^2 + is_{22}^2 - s_{11}s_{33} - s_{13}s_{33} \text{ and }
p_2' = -2s_{12}^2 - 2s_{22}^2 +2is_{11}s_{33} - 2is_{13}s_{33}.\]
{{Consider
\(q_1 = \frac{2p_1' + ip_2'}{4} = is_{12}^2 + s_{11}s_{33} \text{ and } 
q_2 = \frac{2p_1' - ip_2'}{4} = is_{22}^2 - s_{13}s_{33}\). From here, a last simple linear change of variables mapping $s_{11}$ to $-is_{11}$, $s_{13}$ to $is_{13}$, while the rest are fixed, gives the binomial generators $s_{12}^2-s_{11}s_{33}$ and $s_{22}^2-s_{13}s_{33}$.}} 
\end{example}

\subsection{Discussion}
In \Cref{thm:main}, we saw that the dimension of the symmetry Lie algebra of an irreducible projective variety $V$ provides a necessary condition for $V$ to be toric, and \Cref{ex:toric} signals that the symmetry Lie algebra of a variety can be used to provide sufficient conditions that guarantee that a given variety is toric. So we ask:

\begin{question}\label{ques:toric}
Can symmetry Lie algebras detect when $V$ \emph{is} a toric variety? Alternatively, can symmetry Lie algebras detect when there is a linear change of variables under which a prime ideal is toric?
\end{question}
The question has recently found answer \cite{kahle2024efficiently}. We state the main results in the following remark. 
\begin{remark}\label{rem:algorithm}
Follow up work of Kahle and Vill \cite{kahle2024efficiently}, give a positive answer to our question. They do this by searching for a subalgebra in $\g_I$ that is
simultaneously diagonalizable. For this, they first find a Cartan subalgebra and apply Jordan decomposition to it. A coordinate change that diagonalizes this subalgebra
puts $I$ in binomial form if this is possible. All this is implemented in an algorithm. 
They also note that the ideal doesn't need to be prime for deducing whether there is a linear change of variables that turns $I$ into a binomial ideal. In this case, if $I=\langle p_1(x),\ldots, p_k(x)\rangle$, the Lie algebra is the intersection of the Lie algebras $\g_{[I]_{d_i}}$, where $d_i=\deg(p_i(x))$. They provide an efficient algorithm which outputs a linear transformation which turns the ideal in a toric ideal, when possible. 
\end{remark}
With the algorithms, one can now work to fully classify classes of ideals which are toric under a linear transform, including those coming from statistics. 

We end with an example of a non-toric ideal whose symmetry Lie algebra detects that this ideal is toric under a linear transformation.  

\begin{example}
Consider the ideal $I=\langle x_1^2+x_2^2+x_3^2\rangle\subseteq \C[x_1,x_2,x_3]$. This ideal is clearly toric in variables $y_1=x_1$, $y_2=-x_2+ix_3$ and $y_3=x_2+ix_3$, giving $\langle y_1^2-y_2y_3 \rangle\subseteq \mathbb C[y_1,y_2,y_3]$. 
The matrix recording the linear change of variables is \[B=\begin{pmatrix}
    1 & 0 & 0\\ 0 & -1 & 1 \\0 & i & i
\end{pmatrix}.\]
The symmetry Lie algebra of~$I$ is the~$4$-dimensional vector space with basis
\begin{align*}\label{eq:skew}
\scalemath{0.9}{\mathcal{L}=\left\{\begin{pmatrix}
1 & 0 & 0\\
0& 1 & 0\\
0 & 0 & 1
\end{pmatrix}, \begin{pmatrix}
0 & 1 & 0\\
-1& 0 & 0\\
0 & 0 & 0
\end{pmatrix}, \begin{pmatrix}
0 & 0 & 1\\
0& 0 & 0\\
-1& 0 & 0
\end{pmatrix}, \begin{pmatrix}
0 & 0 & 0\\
0& 0 & 1\\
0 & -1 & 0
\end{pmatrix}\right\}.}
\end{align*}
Notice that $B$ changes the basis $\mathcal L$ (so apply  $B^{-1}AB$ to each element in the basis) to the list 
\begin{align*}
\scalemath{0.9}{B^{-1}\mathcal{L}B=
\left\{\left(\begin{array}{rrr}
1 & 0 & 0\\
0& 1 & 0\\
0 & 0 & 1
\end{array}\right), \left(\begin{array}{rrr}
0 & -1 & 1 \\
0.5 & 0 & 0 \\
-0.5 & 0 & 0
\end{array}\right),
\left(\begin{array}{rrr}
0 & i & 1 \\
0.5i & 0 & 0 \\
0.5i & 0 & 0
\end{array}\right),
\left(\begin{array}{rrr}
0 & 0 & 0 \\
0 & -i& 0 \\
0 & 0 & i
\end{array}\right) \right\},}
\end{align*}
which realizes the embedded $2$-dimensional  torus. 
\end{example}
Another question of interest is, given an ideal $I$, to study properties of the ideal associated with the symmetry Lie group $G_I$. As \Cref{ex:lie group manifold} illustrates, they are defined by equations of the same degree. What properties of $I$ are shared with the ideal of $G_I$ and algebra $\g_I$?

\section*{Acknowledgments}
This project initiated at the Texas Algebraic Geometry Symposium 2022 held at Texas A\&M University. The authors are grateful to JM Landsberg for suggesting the problem and for mentoring. We thank Thomas Yahl and Nicklas Day for the useful discussions, and Lisa Nicklasson for comments on a previous version of this paper. We deeply thank the reviewers for the thorough feedback and suggestions on how to make the paper better. 

\bibliographystyle{plain}
\bibliography{ref}

\end{document}

%% file: micro.tex
\usepackage{setspace}
\usepackage[scale=0.85]{geometry}
\usepackage[utf8]{inputenc}
\usepackage{amsmath}
\usepackage{mathrsfs}
\usepackage{amssymb}
\usepackage{graphicx}
\usepackage[colorinlistoftodos]{todonotes}
\usepackage[english]{babel}
\usepackage{amsthm}
\usepackage{titling}
\usepackage{textcomp}
\usepackage{ytableau}
\usepackage{tikz-cd}
\usepackage{eqnarray,amsmath}
\usepackage{xcolor}
\usepackage[colorlinks=true,hyperindex, linkcolor=magenta, pagebackref=false, citecolor=cyan]{hyperref}
\usepackage{graphics}
\usepackage[export]{adjustbox}
\usepackage{subcaption}
\usepackage{epstopdf}
\usepackage{epsfig}
\usepackage{lastpage}
\usepackage{comment}
\usepackage{etoolbox}
\usepackage{mathtools}
\usepackage{cleveref}
\usepackage{ tipa }

\usepackage{url}
\usepackage{microtype}
\usepackage[alphabetic,lite,backrefs]{amsrefs}

\usepackage{amsfonts}
\usepackage{graphicx}
\usepackage{tikz-cd} 
\usetikzlibrary{plotmarks}
\usetikzlibrary{arrows,automata}
\definecolor{darkgreen}{rgb}{0,0.5,0}
\usepackage{graphicx}               
\usepackage{amsmath,amssymb}

\theoremstyle{plain}

\newtheorem{theorem}{Theorem}
\newtheorem{example}[theorem]{Example}
\newtheorem{lemma}[theorem]{Lemma}
\newtheorem{definition}[theorem]{Definition}
\newtheorem{proposition}[theorem]{Proposition}
\newtheorem{corollary}[theorem]{Corollary}
\newtheorem{conjecture}[theorem]{Conjecture}
\newtheorem{remark}[theorem]{Remark}
\newtheorem{question}[theorem]{Question}

\newcommand{\mb}{\mathbb}

\newcommand{\C}{\mathbb C}

\newcommand{\R}{\mb R}

\newcommand{\N}{\mb N}

\newcommand{\GL}{\text{GL}}
\newcommand{\g}{\mathfrak{g}}